\documentclass[fleq]{amsart}
\usepackage{amssymb,amsmath,latexsym,amsbsy,color,enumerate}
\numberwithin{equation}{section}

\def\RR{\mbox{$I\hspace{-.06in}R$}}

\newtheorem{theorem}{Theorem}

\newtheorem{lemma}{Lemma}

\begin{document}
 \title{Approximate solutions to the Dirichlet problem for harmonic maps between hyperbolic spaces}
\author{Duong Minh Duc}
    \address{Department of Mathematics, University of natural sciences, Hochiminh city,
    Vietnam}
 \email{dmduc@hcmuns.edu.vn}
\author{Truong Trung Tuyen}
    \address{Department of Mathematics, Indiana University, Rawles Hall, Bloomington, IN 47405}
 \email{truongt@indiana.edu}
\thanks{This work has been initiated when the second author was at Department of mathematics, University of natural sciences, Hochiminh city, Vietnam.
 He would like to thank Professor Dang Duc Trong for his many invaluable helps. He also would like to express his thankfulness to Professor F. Helein, Professor R. Schoen, and Mr. Le Quang Nam for their generous
 help.}
    \date{\today}
    \keywords{Dirichlet problems; Harmonic functions; Hyperbolic spaces}
    \subjclass[2000]{53A35.}
    \begin{abstract}
Our main result in this paper is the following: Given $H^m,~H^n$ hyperbolic spaces of dimensional $m\geq 2$ and $n$ corresponding, and given a Holder
function $f=(f^1,...,f^{n-1}):\partial H^m\rightarrow
\partial H^n$ between geometric boundaries of $H^m$ and $H^n$. Then for each $\epsilon >0$ there
exists a harmonic map $u:H^m\rightarrow H^n$ which is continuous up
to the boundary (in the sense of Euclidean) and $u|_{\partial
H^m}=(f^1,...,f^{n-1},\epsilon)$.
\end{abstract}
\maketitle
\section{Introduction}
\label{introduction} Let $H^m$ and $H^n$ are hyperbolic spaces with
dimensions $m\geq 2$ and $n$ correspondingly. For convenience, we
use the upper-half space models for $H^m$ and $H^n$. So
$H^m=\{(x^1,...,x^m)\in \RR ^m:~x^m>0\}$, $H^n=\{(y^1,...,y^n)\in
\RR ^n:~y^n>0\}$ with metrics
\begin{eqnarray*}
d_{H^m}^2&=&\frac{1}{(x^m)^2}((dx^1)^2+...+(dx^m)^2),\\
d_{H^n}^2&=&\frac{1}{(y^n)^2}((dy^1)^2+...+(dy^n)^2).
\end{eqnarray*}

So the tension fields of $u=(y^1,...,y^n)$ is
\begin{eqnarray*}
\tau ^{\alpha}=(x^m)^2(\Delta
_0y^{\alpha}-\frac{m-2}{x^m}\frac{\partial y^{\alpha}}{\partial
x^m}-\frac{2}{y^n}<\nabla _0y^{\alpha},\nabla _0y^n>),
\end{eqnarray*}
for $1\leq \alpha\leq n-1$ and
\begin{eqnarray*}
\tau ^n(u)=(x^m)^2(\Delta _0y^n-\frac{m-2}{x^m}\frac{\partial
y^{\alpha}}{\partial x^m}+\frac{1}{y^n}(\sum _{\alpha
=1}^{n-1}|\nabla _0y^{\alpha}|^2-|\nabla _0y^n|^2)),
\end{eqnarray*}
where $\nabla _0$ is the Euclidean gradient and $\Delta _0$ is the
Euclidean Laplacian.

A $C^2$ map $u:~H^m\rightarrow H^n$ is called a harmonic map if
$\tau (u)^s=0$ for all $s=1,2,...,n$. The literature about harmonic
maps between Riemannian manifolds are abundant, we refer the readers
to the classical work \cite{es}.

One of the interesting problems for harmonic maps is that of the
Dirichlet problem at infinity: Given $\partial H^m$ and $\partial
H^n$ geometric boundaries of $H^m$ and $H^n$, and given a continuous
map $f:~\partial M\rightarrow
\partial N$ (here continuity is understood in the sense of
Euclidean), is there a harmonic map $u:~H^m\rightarrow H^n$ such
that in Euclidean sense $u$ is continuous up to the boundary
$\partial H^m$ and takes boundary value $f$?

For this problem with some more requirements for the smoothness of
$f$, there are many results. In three papers \cite{lt1}, \cite{lt2}
and \cite{lt3}, Li and Tam established the existence and uniqueness
of a harmonic function $u$ which is $C^1$ up to the boundary and has
boundary value $f$, provided $f$ is $C^1$. But for more general
types of $f$, according to our knowledge, there is no answer to the
existence of a solution $u$.

In this paper we establish the existence of approximate solutions to the Dirichlet problem for harmonic maps between two hyperbolic spaces with
prescribed boundary value. More explicitly, we prove the following result

\begin{theorem}
Let $f:~H^m\rightarrow H^n$ be a bounded uniformly continuous. Let functions $g$ and $\varphi $ be as in Section 2. Assume that $\int
_0^1t^{-1}g(t)dt<\infty $, in particular, this condition is satisfied if $f$ is Holder continuous. For each $\epsilon
>0$, there exists a harmonic map $u_{\epsilon}:~H^m\rightarrow H^n$
which is continuous up to the boundary $\partial H^m$ and
$u|_{\partial H^m}=(f^1,...,f^{n-1},\epsilon )$.
\label{Holdercasetheorem}\end{theorem}

Our strategy for proving this result is the follows: First, we construct an initial map, i.e., a $C^2$ map $v=(v^1,...,v^{n-1},v^n):H^m\rightarrow H^n$
which has boundary value $f$ for any continuous map $f:\partial H^m\rightarrow H^n$. For this step we follows the ideas in \cite{lt2}, with some changes:
Since the function $f$ needs not to differentiable, we can not take $v^n$ as in \cite{lt2}, and the function $v^n$ of ours is a function of one variable
$x^m$. Then, we use this function to produce harmonic maps $u_{\epsilon}:~H^m\rightarrow H^n$ which takes boundary value $(f^1,...,f^{n-1},v^n+\epsilon
)$ for every $\epsilon >0$.

\section{Initial maps}
\label{initial} In this part, we use the techniques in \cite{lt2} to
construct good initial maps $v$ having the map $f:~\partial
H^m\rightarrow \partial H^n$ as the boundary value.

Let $f:~\RR ^{m-1} \rightarrow \RR ^{n-1}$ be a uniformly continuous
bounded function. Let $g:~H^m\rightarrow (0.\infty )$ be $C^2$,
bounded and
\begin{eqnarray*}
\lim _{x^m\rightarrow 0}g(x',x^m)=0,
\end{eqnarray*}
uniformly in $x'$.

We denote by $v=\{f,g\}:H^m\rightarrow H^n$ the extension of $f$
defined as follows
\begin{eqnarray*}
v^{\alpha}(x',x^m)=\frac{2}{m\omega _m}\int
_{R^{m-1}}\frac{x^mf^{\alpha}(y')}{(|x'-y'|^2+(x^m)^2)^{m/2}}dy',
\end{eqnarray*}
for $1\leq \alpha \leq n-1$ and
\begin{eqnarray*}
v^{n}(x',x^m)=g(x',x^m).
\end{eqnarray*}

By results in \cite{lt2} (pp. 628-630) we have

(i) $v$ is $C^2$ and up to the boundary given by $x^m=0$ it is
continuous.

(ii) If $1\leq \alpha \leq n-1$ then
\begin{eqnarray*}
\lim _{x^m\rightarrow 0}x^m|\nabla _0v^{\alpha}|=0,
\end{eqnarray*}
uniformly in $x'$.

Moreover, by estimates of elliptic PDEs (see Theorem 2.10 in
\cite{gt}), noting that $v^{\alpha}$ is bounded, there exists
constants $C>0$ such that
\begin{equation}
\max\{(x^m)^3|D^{3}v^{\alpha}|,(x^m)^2|D^2 v^{\alpha}|,(x^m)|\nabla
_0v^{\alpha}|\}\leq C. \label{importantproperty}
\end{equation}

We put
\begin{eqnarray*}
g(r)=\sup _{x',y'\in\RR ^{m-1},~|x'-y'|\leq r}|f(y')-f(x')|,
\end{eqnarray*}
and
\begin{eqnarray*}
\varphi (r)=\int _{0}^{\infty}\frac{r}{s^2+r^2}g(s)ds.
\end{eqnarray*}
Since $g$ is monotone it follows that $g$ is Lebesgue measurable.
Moreover, since $g$ is bounded, we see that $\varphi $ is
well-defined.

Using polar coordinates with center at $x'$ we see that there exists
a constant $C>0$ such that
\begin{eqnarray*}
\int _{\RR
^{m-1}}\frac{x^m|f(y')-f(x')|}{(|x'-y'|^2+(x^m)^2)^{m/2}}\leq
C\varphi (x^m),
\end{eqnarray*}
for all $x'\in \RR ^{m-1}$.

Since $f$ is uniformly continuous we see that
\begin{eqnarray*}
\lim _{r\rightarrow 0}g(r)=0.
\end{eqnarray*}

Now we show that
\begin{eqnarray*}
\lim _{x^m\rightarrow 0}\varphi (x^m)=0.
\end{eqnarray*}

Indeed, for any $\epsilon >0$, we find $\delta >0$ such that
\begin{eqnarray*}
g(s)\leq \epsilon,
\end{eqnarray*}
if $0<s\leq \delta$. So, if $K=\sup _{s\in\RR}g(s)$ we have
\begin{eqnarray*}
\varphi (r)&=&\int _{0}^{\delta}\frac{r}{s^2+r^2}g(s)ds+\int _{\delta}^{\infty}\frac{r}{s^2+r^2}g(s)ds\\
&\leq&\int _{0}^{\delta}\epsilon\frac{r}{s^2+r^2}ds+\int _{\delta}^{\infty}K\frac{r}{s^2+r^2}ds\\
&=&\epsilon \arctan (\delta /r)+K(\pi/2 -\arctan (\delta /r)).
\end{eqnarray*}
Letting $r\rightarrow 0$ we see that
\begin{eqnarray*}
\limsup _{r\rightarrow 0}\varphi (r)\leq \epsilon \pi /2.
\end{eqnarray*}
Since $\epsilon >0$ is arbitrary, we see that
\begin{eqnarray*}
\lim _{r\rightarrow 0}\varphi (r)=0.
\end{eqnarray*}

Thus, if we put $v=\{f,\varphi (x^m)\}$ we see that $v$ is an extension of $f$. Moreover we have the following result

\begin{lemma} Let $f:~\partial H^m\rightarrow \partial H^n$ be nonconstant, uniformly continous and bounded. Put $v=\{f,\varphi (x^m)\}$ as above. Then $v$ is smooth, up to the boundary it is continuous, $v|_{\RR ^{m-1}}=f$ and there exists $C>0$ such that
for $x^m$ near $0$ we have
\begin{eqnarray*}
||\tau (v)||^2\leq C.
\end{eqnarray*}
\label{lemmainitialmaps}\end{lemma}
\begin{proof} By Section 6 in \cite{lt2} we have
\begin{eqnarray*}
|(x^m)\nabla _0v^{\alpha}|\leq C_3|\varphi (x^m)|,
\end{eqnarray*}
where $1\leq \alpha\leq n-1$ and $C_3$ is a positive constant.

Directly computation gives
\begin{eqnarray*}
\varphi '(r)&=&\int _0^{\infty}\frac{s^2-r^2}{(s^2+r^2)^2}g(s)ds,\\
\varphi "(r)&=&\int _0^{\infty}\frac{-2r}{(s^2+r^2)^2}g(s)ds+\int
_{0}^{\infty}\frac{-4r(s^2-r^2)}{(s^2+r^2)^3}g(s)ds.
\end{eqnarray*}

So
\begin{eqnarray*}
\max\{|r\varphi '(r)|,|r^2\varphi "(r)|\}\leq C_4\varphi (r),
\end{eqnarray*}
where $C_4$ is a constant.

Since $g$ is increasing, $g'$ exists almost everywhere and $g'\geq
0$. Using integration by parts, noting that
$\frac{d}{ds}(\frac{-s}{s^2+r^2})=\frac{s^2-r^2}{(s^2+r^2)^2}$, we
have
\begin{eqnarray*}
\varphi '(r)&=&\int _0^{\infty}\frac{s^2-r^2}{(s^2+r^2)^2}g(s)ds\\
&=&\frac{-s}{s^2+r^2}g(s)|_0^{\infty}+\int _0^{\infty}\frac{s}{s^2+r^2}g'(s)ds\\
&=&\int _0^{\infty}\frac{s}{s^2+r^2}g'(s)ds.
\end{eqnarray*}

Differentiating the last term in above equality we get
\begin{eqnarray*}
\varphi "(r)=-\int _0^{\infty}\frac{2rs}{(s^2+r^2)^2}g'(s)ds.
\end{eqnarray*}

Since $f$ is nonconstant we see easily that $g'\not\equiv 0$ (in
fact, we don't need this restriction since we can add $g$ with a
non-constant positive function, for example $(x^m)^{1/2}$ ). So
since $g'\geq 0$, it follows from above equalities that
\begin{eqnarray*}
\varphi '(r)>0,
\end{eqnarray*}
and
\begin{eqnarray*}
|r\varphi "(r)|\leq C_5\varphi '(r),
\end{eqnarray*}
where $C_5$ is a positive constant. Then use the formula for the tension field we are done.
\end{proof}
\section{Proof of Theorem \ref{Holdercasetheorem}}
\begin{proof}
Fixed $\epsilon >0$. We define $v_{\epsilon }:~H^m\rightarrow H^n$
as follows:
\begin{eqnarray*}
v_{\epsilon}(x)=(v^1(x),v^2(x),...,v^{n-1}(x),\varphi (x^m)+\epsilon
).
\end{eqnarray*}
For each $\delta >0$ denote $u_{\epsilon
,\delta}:~H^m\supseteqq\Omega _{\delta}=\{x^m>\delta \} \rightarrow
H^n$ the harmonic map taking value $v_{\epsilon}$ on $\partial
\Omega _{\delta}$.

By inequality (2.1) in \cite{dw} and properties of $v$ and
$v_{\epsilon}$ (see Lemma \ref{lemmainitialmaps}) we have
\begin{eqnarray*}
\Delta _{H^m}d_{H^n}(u_{\epsilon ,\delta},v_{\epsilon})\geq -|\tau
(v_{\epsilon})|\geq -C\frac{\varphi (x^m)}{\varphi (x^m)+\epsilon
}\geq -C\frac{1}{\epsilon }\varphi (x^m),
\end{eqnarray*}
for all $x\in \Omega _{\delta}$, and here $C$ is one constant from
Lemma \ref{lemmainitialmaps}.

We claim that the function
\begin{eqnarray*}
\psi (r)=\int _{0}^r\int _{s}^{\infty}u^{-2}\varphi (u)~du~ds
\end{eqnarray*}
is well-defined for $r\geq 0$. In fact, using the formula for
$\varphi$ we have
\begin{eqnarray*}
\psi (r)=\int _{0}^r\int _{s}^{\infty}u^{-2}\varphi ~du~ds=\int
_0^r\int _s^{\infty}\int
_{0}^{\infty}u^{-1}(u^2+t^2)^{-1}g(t)dt~du~ds.
\end{eqnarray*}
Since the integrand is non-negative, using Fubini's theorem we have
\begin{eqnarray*}
\int _0^r\int _s^{\infty}\int
_{0}^{\infty}u^{-1}(u^2+t^2)^{-1}g(t)dt~du~ds&=&\int _0^r\int
_{0}^{\infty}\int _s^{\infty}u^{-1}(u^2+t^2)^{-1}g(t)du~dt~ds\\
&=&\int _0^r\int _{0}^{\infty}\frac{1}{2t^2}\log
(1+\frac{t^2}{s^2})g(t)dt~ds\\
&=&\int _{0}^{\infty}\int _0^r\frac{1}{2t^2}\log
(1+\frac{t^2}{s^2})g(t)ds~dt\\
&=&\int _{0}^{\infty}\frac{\pi t^2 -2\arctan (\frac{t}{r})t^2+r\log
(1+\frac{t^2}{r^2})}{t^3}g(t)dt.
\end{eqnarray*}
Now since $g(t)$ is bounded we have
\begin{eqnarray*}
\int _{0}^{\infty}\frac{\pi t^2 -2\arctan
(\frac{t}{r})t^2}{t^3}g(t)dt
\end{eqnarray*}
is convergent. Fixed $r\geq 0$, near $t=0$ we have
\begin{eqnarray*}
\frac{r\log (1+\frac{t^2}{r^2})}{t^3}g(t)\approx t^{-1}g(t),
\end{eqnarray*}
and when $t\rightarrow\infty$ we have
\begin{eqnarray*}
\frac{r\log (1+\frac{t^2}{r^2})}{t^3}g(t)\approx t^{-3}g(t),
\end{eqnarray*}
hence since $g(t)$ is bounded and the assumption that $\int
_0^1t^{-1}g(t)$ converges, our claim is verified.

We use the same $\psi $ to denote the function $\psi :~H^m\rightarrow \RR$ defined by $\psi (x)=\psi (x^m)$ for $x=(x^1,...,x^{m-1},x^m)\in H^m$. Now we
have $\psi '(r)=\int _{r}^{\infty}u^{-2}\varphi ~du
>0$ and $\psi "(r)=-r^{-2}\varphi (r)$, since $m\geq 2$ we have
\begin{eqnarray*}
\Delta _{H^m}(-\psi (x))=-(x^m)^{2}[\psi "(x^m)-\frac{(m-2)}{x^m}\psi '(x^m)]\geq -(x^m)^{2}\psi "(x^m)=\varphi (r).
\end{eqnarray*}
Hence
\begin{eqnarray*}
\Delta _{H^m}(d_{H^n}(u_{\epsilon
,\delta},v_{\epsilon})-C\frac{1}{\varphi (\epsilon )}\psi )\geq 0,
\end{eqnarray*}
for $x\in \Omega _{\delta}$. Hence by maximum principle we have
\begin{eqnarray*}
\sup _{x\in \Omega}d_{H^n}(u_{\epsilon ,\delta},v_{\epsilon})\leq
C\frac{1}{\epsilon }\psi (x^m).
\end{eqnarray*}
This bound for $d_{H^n}(u_{\epsilon ,\delta},v_{\epsilon})$ is
independent of $\delta$, hence by standard arguments (see the proof
of Theorem 6.4 in \cite{lt2}) we have a harmonic map
$u_{\epsilon}:~H^m\rightarrow H^n$ which is the subsequent limit of
$u_{\epsilon ,\delta}$. Moreover for all $x\in H^m$ we have
\begin{eqnarray*}
d_{H^n}(u_{\epsilon},v_{\epsilon})\leq C\frac{1}{\epsilon }\psi
(x^m).
\end{eqnarray*}
Hence
\begin{eqnarray*}
\lim _{x^m\rightarrow 0}d_{H^n}(u_{\epsilon},v_{\epsilon})=0,
\end{eqnarray*}
which shows that $u_{\epsilon }$ is continuous up to the boundary
and takes boundary value $v_{\epsilon
}(x^1,...,x^{m-1},0)=(f^1,f^2,...,f^{n-1},\epsilon )$.
\end{proof}

\end{document}